\documentclass[leqno,12pt]{amsart} 
\usepackage{amssymb, enumitem, xcolor}
\newcommand{\R}{\mathbb{R}}
\newcommand{\N}{\mathbb{N}}
\theoremstyle{plain} 
\newtheorem{theorem}{\indent\sc Theorem}[section]
\newtheorem{lemma}[theorem]{\indent\sc Lemma}

\theoremstyle{definition} 

\begin{document}

\title[FDE: uniqueness of solutions with a moving singularity]{Fast diffusion equation: uniqueness of solutions with a moving singularity}

\author[M. Fila]{Marek Fila \textdagger} 

\author[P. Mackov\'a]{Petra Mackov\'a} 


%
\thanks{ 
The authors would like to thank Jin Takahashi for his valuable comments and discussion. 
The first author was partially supported by the Slovak Research and Development Agency under contract No. APVV-18-0308 and by  VEGA grant 1/0339/21. 
The second author was partially supported by VEGA grant 1/0339/21 and Comenius University grant UK/236/2022 and UK/242/2023. 
\textdagger 
Sadly, Marek Fila, a supervisor, friend, and co-author of this paper, passed away in April 2023. 
In his memory, the second author has decided to publish the research, as this work with his valuable impact was finished before his passing.}
\address{
Department of Applied Mathematics and Statistics \endgraf
Comenius University \endgraf
842 48 Bratislava \endgraf
Slovakia
}

\address{
Department of Applied Mathematics and Statistics \endgraf
Comenius University \endgraf
842 48 Bratislava \endgraf
Slovakia
}
\email{petra.mackova@fmph.uniba.sk}

\maketitle

\begin{abstract}

We focus on open questions regarding the uniqueness of distributional solutions of the fast diffusion equation (FDE) with a given source term. 
When the source is sufficiently smooth, the uniqueness follows from standard results.
Assuming that the source term is a measure, the existence of different classes of solutions is known, but in many cases, their uniqueness is an open problem.
In our work, we focus on the supercritical FDE and prove the uniqueness of distributional solutions with a Dirac source term that moves along a prescribed curve. 

\end{abstract}

\section{Introduction}

Let $0<m<1$, $n \geq 3$, and $0<T\leq \infty$. We study the uniqueness of distributional solutions of the fast diffusion equation
\begin{equation} \label{eq:nonlinear}
    u_t = \Delta u^m + f(x,t), \quad \text{ in } \quad \mathcal{D}'(\R^n \times (0, T)),
\end{equation}
where $f$ is a given source term.
More specifically, we are interested in solutions of~\eqref{eq:nonlinear} that satisfy $u \in L^1_{loc}(\R^n \times (0,T) )$ and the integral equality
\begin{equation} \label{eq:distrib}
    \int_0^T \int_{\R^n} \big( u \varphi_t + u^m \Delta\varphi + f \varphi\big) \, dx \, dt =  0
\end{equation}
for all $\varphi \in C_0^\infty (\R^n \times (0, T))$.
If, moreover, $\nabla u^m \in L^1_{loc}(\R^n \times (0,T) )$ and $u$ satisfies
\[
    \int_0^T \int_{\R^n} \big( u \varphi_t - \nabla u^m \cdot \nabla\varphi + f \varphi\big) \, dx \, dt = 0
\]
for all $\varphi \in C_0^\infty (\R^n \times (0, T))$ then we call it a weak solution of~\eqref{eq:nonlinear}.

Some techniques to prove uniqueness of solutions of~\eqref{eq:nonlinear} can be found in the book~\cite{vazquez1} by V{\'a}zquez.
Focusing on weak solutions and assuming that $u \in L_{loc}^2(\R^n \times (0,T))$, $u^m \in L_{loc}^2(0,T;H_0^1(\R^n))$, and $f \in L_{loc}^1(\R^n \times (0,T))$, one can use a test function $\varphi$ of the form
\begin{equation*}
	\varphi(x,t) = \begin{cases} 
		\int_t^T \big( u_1^m(x,s)-u_2^m(x,s)\big) ds & \text{if } 0<t<T, \\
		0 & \text{if } t \geq T,
	\end{cases}
\end{equation*}
which was introduced by Oleinik~\cite{Oleinik1957}.



The critical exponent $m_c:=(n-2)/n$ plays an important role in the theory of the fast diffusion equation. See, for example V{\'a}zquez~\cite{vazquez2}.
In~\cite{Lu2}, Lukkari studies solutions of the fast diffusion equation in the range $m_c<m<1$ with $\Omega$ instead of $\R^n$, where $\Omega$ is a bounded domain with a smooth boundary. Assuming that the forcing term $f$ is a non-negative Radon measure on $\R^{n+1}$ such that $f(\Omega \times (0,T)) < \infty$, he proves the existence of a specific class of weak solutions of~\eqref{eq:nonlinear} in cylinders of the form $\Omega \times (0,T)$. These solutions satisfy $u \in L^q((0,T); W_0^{1,q}(\Omega))$, where $q$ is any number such that $1 \leq q < 1 + 1/(1 + mn)$. Since the upper bound on $q$ is always less than $2$, Lukkari's weak solutions lack the $L^2$-integrability conditions assumed by V{\'a}zquez in~\cite{vazquez1}, hence, their uniqueness was left as an open problem.

A standard uniqueness result for $0<m<1$ by Herrero and Pierre can be found in~\cite{HP}. Here, the authors prove the uniqueness of distributional solutions of the signed fast diffusion equation, i.e. 
\begin{equation} \label{eq:distr_signed}
 u_t = \Delta (u|u|^{m-1}) \quad \text{ in } \quad \mathcal{D}'(\R^n \times (0, \infty)),
\end{equation}
assuming that $u\in C([0,\infty);L^1_{loc}(\R^n))$ and the time derivative satisfies $u_t\in L^1_{loc}(\R^n \times (0,\infty))$.

More recently, new results concerning uniqueness of subcritical fast diffusion have been found. 
In~\cite{TY21}, Takahashi and Yamamoto focused on the case when $n \geq 3$ and $0<m<m_c$. 
They showed the uniqueness of signed solutions of the initial value problem 
\begin{align}
	\label{eq:main_standing}
	u_t &= \Delta (u|u|^{m-1}), \qquad x \in \R^n \setminus \{ \xi_0 \}, \quad t \in (0, T), \\
	\label{eq:main_initial_standing}
	u(x,0)&=u_0(x), \qquad \qquad \, \, \, x \in \R^n \setminus \{\xi_0\},
\end{align}
with $0<T\leq \infty$ and $\xi_0 \in \R^n$.
More precisely, the authors of~\cite{TY21} proved that for two functions $u_1, u_2$ such that $u_1, u_2 \in C^{2,1}((\R^{n}\setminus\{\xi_0\})\times(0,T)) \cap C( (\R^{n}\setminus\{\xi_0\})\times[0,T))$ that satisfy~\eqref{eq:main_standing}-\eqref{eq:main_initial_standing} pointwise
and $u_1(\cdot,0)=u_2(\cdot,0)$ on $\R^{n}\setminus\{\xi_0\}$,
it holds that $u_1 \equiv u_2$ on $(\R^{n}\setminus\{\xi_0\})\times(0,T)$.
Hui demonstrated in~\cite{Hui2020} that if $n \geq 3$ and $0<m<m_c$, under suitable conditions on initial data, solutions that have a finite number of standing singularities are also uniquely determined. By solutions with finitely many standing singularities, we mean that these solutions satisfy equations~\eqref{eq:main_standing}-\eqref{eq:main_initial_standing} with $\R^n \setminus \{ \xi_0, \xi_1, \dots, \xi_{i} \}$ instead of $\R^n \setminus \{ \xi_0 \}$. Here, $i\in\N$ and $\xi_0, \xi_1, \dots, \xi_{i} \in \R^n$.

More generally, we can assume that $\xi: [0,T) \rightarrow \R^n$ is a given curve and study the problem
\begin{align}
	\label{eq:main}
	u_t &= \Delta (u|u|^{m-1}), \qquad x \in \R^n \setminus \{ \xi(t) \}, \quad t \in (0, T), \\
	\label{eq:main_initial}
	u(x,0)&=u_0(x), \qquad \qquad \, \, \, x \in \R^n \setminus \{\xi(0)\},
\end{align}
with a moving singularity $\xi(t)\not\equiv \xi(0)$ for some $t\in(0,T)$ and $0<T\leq \infty$. %
In the case $m>m_c$ and $T=\infty$, positive asymptotically radially symmetric solutions of the initial value problem~\eqref{eq:main}-\eqref{eq:main_initial} were studied in~\cite{FTY, FMTY, FMTY21}. 
Along the given curve $\xi$ with suitable properties, these solutions keep a singularity at all times, i.e. 
$u(x,t) \to \infty$ as $x \to \xi(t)$ for each $t\in[0,T)$.
Our main result concerns the uniqueness of these solutions in the supercritical fast diffusion case. In order to formulate this result, we give a precise description of solutions from~\cite{FTY, FMTY21}. 
Let $n \geq 3$ and $T=\infty$. Assume that either
\begin{enumerate} [label={(A\upshape\arabic*)}, start=0, align=left, widest=iii, leftmargin=*]
    \item \label{eq:xi0} $m > m_c$ and $\xi(t) \equiv \xi_0$ for some $\xi_0 \in \R^n$,
\end{enumerate}
or
\begin{enumerate} [label={(A\upshape\arabic*)}, align=left, widest=iii, leftmargin=*]
    \item \label{eq:xi} $m > m_*:= (n-2)/(n-1)$ and $\xi \in C^1([0,\infty);\R^n)$, $\xi'$ is locally H\"older continuous, and there exist positive constants $\Xi, \beta$ such that $|\xi'(t)|\leq \Xi e^{-\beta t}$ for $t \geq 0$.
\end{enumerate}
Assume, moreover, that
\begin{enumerate} [label={(A\upshape\arabic*)}, align=left, start=2, widest=iii, leftmargin=*]
    \item \label{eq:k} $k\in C^1([0,\infty))$ satisfies $\kappa^{-1}\leq k(t)\leq \kappa$ and $|k'(t)|\leq \kappa'$ for $t \geq 0$ and some positive constants $\kappa$ and $\kappa'$,
    \item \label{eq:u0} $u_0(x) \in C(\R^n \setminus \{\xi(0) \} )$ is positive and there exist $\lambda$, $\mu$ and $\nu$ satisfying
\begin{equation} \label{eq:const}
    \max \left\{ (n-2)/m -1, 0 \right\} < \lambda <\mu< n-2<\nu
\end{equation}
    such that $u_0(x)^m = k(0)^m |x-\xi(0)|^{-n+2} + O(|x-\xi(0)|^{-\lambda})$ as $x\to\xi(0)$, and $C^{-1} |x-\xi(0)|^{-\nu} \leq u_0(x)^m \leq C |x-\xi(0)|^{-\mu}$ for $|x-\xi(0)|\geq 1$ with some constant $C>1$.    
\end{enumerate}

Under these assumptions, \cite{FTY} implies the existence of a function $u>0$ satisfying the following:
\begin{enumerate} [label={\upshape(\roman*)}, align=left, widest=iii, leftmargin=*]
    \item \label{c:1} $u\in C^{2,1}( \{(x,t)\in \R^{n+1}: x\neq \xi(t), \, t \in (0,\infty)\}) \cap C( \{(x,t)\in \R^{n+1}: x\neq \xi(t), \, t \in [0,\infty)\})$ and $u>0$ satisfies \eqref{eq:main}-\eqref{eq:main_initial} pointwise,
    \item \label{c:l1}  $u\in C([0,\infty);L^1_{loc}(\R^n))$,
    \item \label{c:asymp} for each $t\geq 0$, $u$ has the asymptotic behavior
    \begin{equation*} 
	u(x,t)^m = k(t)^m |x-\xi(t)|^{-n+2} + O(|x-\xi(t)|^{-\lambda}) \quad \text{ as } \quad x\to\xi(t),
    \end{equation*}
    \item \label{c:4} for $t \geq 0$ and $|x-\xi(t)|\geq 1$, it holds that
    \begin{equation*}
	C^{-1} e^{-Ct} |x-\xi(t)|^{-\nu} \leq u(x,t)^m \leq C e^{Ct} |x-\xi(t)|^{-\mu}
    \end{equation*}
     with some constant $C>1$.
\end{enumerate}
We note that \cite{FTY} dealt with moving singularities, i.e. the existence was proved under assumptions~\ref{eq:xi}, \ref{eq:k}, \ref{eq:u0}. Later, in~\cite{FMTY21} it was remarked that the existence from~\cite{FTY} is valid in the whole supercritical parameter range $m>m_c$ if the singularity is standing (i.e. assuming~\ref{eq:xi0}, \ref{eq:k}, \ref{eq:u0}). 

Moreover, it was established in~\cite{FMTY21} that a function $u$ from~\cite{FTY} satisfying~\ref{c:1}-\ref{c:asymp} is a distributional solution of problem~\eqref{eq:nonlinear} with a weighted moving Dirac source term
\begin{equation*} \label{eq:n3_weak_delta}
    u_t = \Delta u^m + (n-2) |S^{n-1}| k^m(t) \delta_{\xi(t)}(x) \quad \text{ in } \quad \mathcal{D}'(\R^n \times (0, \infty)).
\end{equation*}
More precisely, $u$ satisfies~\eqref{eq:distrib} for all $\varphi \in C_0^\infty (\R^n \times (0, \infty))$ with $$f(x,t)=(n-2) |S^{n-1}| k^m(t) \delta_{\xi(t)}(x).$$
Here, $\delta_{\xi(t)}$ gives unit mass to the point $\xi(t) \in \R^n$ for each $t\geq 0$, and $|S^{n-1}|$ denotes the surface area of the $(n-1)$-dimensional unit sphere. 
A Dirac measure that moves with time can be also found as a source in parabolic systems, and this phenomenon has been used to model various biological scenarios, such as axon growth or angiogenesis, as discussed in~\cite{CZ} and~\cite{Bookholt}, respectively.
We summarize our main result in the theorem below.

\begin{theorem} \label{th:uniqueness}
    Let $n \geq 3$ and $T=\infty$. 
    Assume that either~\ref{eq:xi0} or~\ref{eq:xi} holds. Assume, moreover, that conditions~\ref{eq:k}, and~\ref{eq:u0} are satisfied, and that functions $u_1, u_2$ satisfy~\ref{c:1}--\ref{c:asymp}. 
    Then the equality $u_1(\cdot,0)=u_2(\cdot,0)$ on $\R^{n}\setminus\{\xi(0)\}$ implies that $u_1 \equiv u_2$ on $\{(x,t)\in \R^{n+1}: x\neq \xi(t), \, t \in (0,\infty) \}$.
\end{theorem}

Uniqueness results for the porous medium equation can be found in~\cite{DK, P, vazquez1}.
For the uniqueness of solutions of a semi-linear parabolic equation with singularity moving along a prescribed curve, see~\cite{SY2009}, where similar conditions to~\ref{c:1}-\ref{c:4} were considered. For non-uniqueness examples for a semilinear heat equation, see e.g. \cite{FMY2016} and references therein.

The remainder of this paper is dedicated to the proof of Theorem~\ref{th:uniqueness}.

\section{Proof of Theorem~\ref{th:uniqueness} } 
\label{sec:uniqueness_m<1}

\begin{proof}
This proof is based on ideas of Takahashi, Yamamoto, F., M., Yanagida, Herrero, and Pierre, see \cite{TY21}, \cite{FMTY21}, and~\cite{HP}.

\textit{Step 1.}
Set sign$(f)=f/|f|$ for $f \neq 0$ and sign$(f)=0$ for $f = 0$.
We recall that for a locally integrable function $f$ such that $\Delta f \in L^1_{loc}(D)$ in $D\subseteq \R^n$, Kato proved the distributional inequality
\begin{equation*} \label{eq:kato}
    \text{sign}(f) \Delta f \leq \Delta|f|.
\end{equation*}
Let $u_1, u_2$ be two functions satisfying assumptions~\ref{c:1}-\ref{c:asymp} and $u_1(\cdot,0)=u_2(\cdot,0)$ on $\R^{n}\setminus\{\xi(0)\}$.
Then it holds that
\begin{equation} \label{eq:using_kato}
    \partial_t |u_1-u_2| = \text{sign}(u_1-u_2) \partial_t (u_1-u_2) = \text{sign}(u_1-u_2) \Delta (u_1^m-u_2^m) \leq \Delta |u_1^m-u_2^m|
\end{equation}
for $x \in \R^n \setminus \{ \xi(\tau) \}$ and $\tau \in [0, \infty)$. 

The choice of test function is based on Lemma 2.1 from~\cite{TY15}. We present this lemma for completeness.
\begin{lemma} \label{lemma}
    \cite{TY15}. Let $n\geq 1$, $\tau_1,\tau_2 \in \R$, $\tau_1<\tau_2$, and $\alpha \in (0,1]$. Suppose that $\xi(\tau)$ is locally H\"older continuous in $\tau\in\R$ for some  $\alpha \in (0,1]$. Then there exist $\varepsilon_0~=~\varepsilon_0(n,\alpha,\tau_1,\tau_2) \in (0,1)$ and $c_0=c_0(n,\alpha,\tau_1,\tau_2)>0$ independent of $x,\tau, \varepsilon$ with the following property: for any $\varepsilon\in (0,\varepsilon_0)$ there exists a family of cut-off functions $\{ \eta_\varepsilon \}_{\varepsilon>0} \subset C^\infty (\R^n \times \R)$ such that $0\leq \eta_\varepsilon \leq 1$,
    \begin{equation*}
	\eta_\varepsilon= \begin{cases} 
		0 & \text{if } |x-\xi(\tau)|< \varepsilon, \\
		1 & \text{if } |x-\xi(\tau)|> 2\varepsilon,
	\end{cases}
    \end{equation*}
    and for $(x,\tau)\in \R^n \times [\tau_1,\tau_2]$,
    \begin{equation*}
        |\nabla \eta_\varepsilon| \leq c_0 \varepsilon^{-1}, \quad 
        |\Delta \eta_\varepsilon| \leq c_0 \varepsilon^{-2}, \quad 
        |\partial_\tau  \eta_\varepsilon| \leq c_0 \varepsilon^{-1/\alpha} \quad 
        \text{on} \quad \R^n \times [\tau_1,\tau_2].
    \end{equation*}
\end{lemma}
We note that since $\xi \in C^1([0,\infty);\R^n)$ by~\ref{eq:xi}, we can extend it to the whole space so that the local H\"older continuity assumption in $\tau\in\R$ in Lemma~\ref{lemma} is satisfied.  
Let $\varphi \in C_0^\infty (\R^n)$ be a nonnegative function and set $\varphi_\varepsilon(x,\tau) := \eta_\varepsilon(x,\tau) \varphi(x)$.
For $R>0$ and $z\in\R^n$, we let $B_R \left(z\right):= \{ x\in \R^n; |x-z|<R \}$.
For simplicity, by $B_R:=B_R \left(\xi(\tau)\right)$ we will denote an open ball with radius $R$ centered at $\xi(\tau)$. We note that $\varphi_\varepsilon=\varphi$ for $x\in \R^n \setminus B_{2\varepsilon}$.

\eqref{eq:using_kato} with $\varphi_\varepsilon$ and integrating both sides by parts, we have
\begin{equation*}
	\partial_\tau \int_{\mathbb{R}^n} \varphi_\varepsilon\,|u_1-u_2|\,dx - \int_{\mathbb{R}^n} (\partial_\tau\eta_\varepsilon)\, \varphi\,|u_1-u_2|\,dx
	\leq 
    \int_{\R^n} |u_1^m-u_2^m| \Delta \varphi_\varepsilon \, dx.
\end{equation*}
Since $\partial_\tau \eta_\varepsilon$ vanishes outside the region
$B_{2\varepsilon}\setminus B_{\varepsilon}$, this gives us
\begin{equation*} \label{eq:proof_beginning}
    \partial_\tau \int_{\R^n} \varphi_\varepsilon |u_1-u_2| \, dx 
    \leq \int_{\R^n} |u_1^m-u_2^m| \Delta \varphi_\varepsilon \, dx 
    + L_\varepsilon,
\end{equation*}
where we denote 
\begin{equation*} \label{eq:est_f_L}
    L_\varepsilon (\tau) :=
    \int_{B_{2\varepsilon}\setminus B_{\varepsilon}}
    \partial_\tau \eta_\varepsilon(x,\tau) \,
    \varphi(x)\,\lvert u_1(x,\tau)-u_2(x,\tau)\rvert\,dx.
\end{equation*}
We fix $t>0$. Since $u_1(\cdot,0)=u_2(\cdot,0)$ on $\R^{n}\setminus\{\xi(0)\}$, 
integrating the above inequality
with respect to $\tau$ from $0$ to $t$ gives
\begin{equation*}
\begin{split}
    \int_{\R^n } \varphi_\varepsilon(x,t) |u_1(x,t)-u_2(x,t)| \, dx 
    \leq &\int_0^t \int_{\R^n} |u_1(x,\tau)^m-u_2(x,\tau)^m| \Delta \varphi_\varepsilon(x,\tau) \, dx \, d\tau \\
    &+ \int_0^t L_\varepsilon (\tau) \, d\tau.
\end{split}
\end{equation*}
This can be written as
\begin{equation*} \label{eq:p6-1-2}
\begin{split}
    \int_{\R^n \setminus B_{2\varepsilon}} \varphi |u_1-u_2| \, dx 
    \leq & \int_0^t \int_{\R^n \setminus B_{2\varepsilon}} |u_1^m-u_2^m| \Delta \varphi \, dx \, d\tau \\
    &+ H_\varepsilon + \int_0^t \left( I_\varepsilon + J_\varepsilon + K_\varepsilon + L_\varepsilon \right)\, d\tau , 
\end{split}
\end{equation*}
where we use similar notation as in~\cite{FMTY21}, i.e. we denote
\begin{align*} \label{eq:IJKeps}
\begin{split}
    H_\varepsilon &:= -\int_{B_{2\varepsilon}\setminus B_\varepsilon} \eta_\varepsilon(x,t) \varphi(x) |u_1(x,t)-u_2(x,t)| \, dx, \\
    I_\varepsilon &:=  \int_{B_{2\varepsilon}\setminus B_{\varepsilon}} |u_1^m-u_2^m| \eta_\varepsilon\Delta\varphi \, dx, \\
    J_\varepsilon &:= 2\int_{B_{2\varepsilon}\setminus B_{\varepsilon}} |u_1^m-u_2^m| \nabla\eta_\varepsilon \cdot \nabla\varphi \, dx,  \\
    K_\varepsilon &:= \int_{B_{2\varepsilon}\setminus B_{\varepsilon}} |u_1^m-u_2^m| \varphi\Delta\eta_\varepsilon \, dx.
\end{split}
\end{align*}

\textit{Step 2.}
In what follows, we make use of the characterization of the behavior of $u_1, u_2$ in a neighborhood of the moving singularity \(\xi(\tau)\).
We want to pass to the limit as $\varepsilon \to 0$ and prove that 
\begin{equation*} \label{eq:assertion}
    H_\varepsilon, \, \, I_\varepsilon, \, \, J_\varepsilon, \, \, K_\varepsilon, \, \, L_\varepsilon \, 
    \to 0
    \quad \text{as } \quad \varepsilon \to 0,
\end{equation*}
uniformly in $\tau \in (0,t)$.
As in~\cite{FMTY21}, we choose $\varepsilon$ sufficiently small so that the method of sub- and supersolutions in~\cite{FTY} provides estimates of the form
\begin{equation} \label{eq:v_estimates}
\begin{split}
    u^m(x,\tau) &\leq k^m(\tau) \left(|x-\xi(\tau)|^{2-n} + b(\tau)|x-\xi(\tau)|^{-\lambda}\right),\\
    u^m(x,\tau) &\geq k^m(\tau) \left(|x-\xi(\tau)|^{2-n} - b(\tau)|x-\xi(\tau)|^{-\lambda}\right)_+,
\end{split}
\end{equation}
for all $x \in B_{2\varepsilon}$ and $\tau \in [0,t]$.
Here, $b(\tau)=b_0e^{B \tau}$ for some constants $B$, $b_0 >1$, $\lambda < n-2$ by~\eqref{eq:const}, and we recall that $k$ is a given function satisfying~\ref{eq:k}.
In what follows, by $c$ we will denote a large enough but otherwise arbitrary constant independent of $t, \tau$ and $\varepsilon$.
Inspecting the proof of Theorem~1.5 in~\cite{FMTY21}, we see that for $\tau\in [0,t]$ we have
\begin{equation*} 
\begin{split}
    |I_\varepsilon| &\leq \int_{B_{2\varepsilon}\setminus B_{\varepsilon}} (u_1^m+u_2^m) \eta_\varepsilon |\Delta\varphi| \, dx 
    \leq c \int_{\varepsilon}^{2\varepsilon} r \, dr \to 0 \quad \text{ as } \quad \varepsilon \to 0, \\
    |J_\varepsilon| &\leq 2\int_{B_{2\varepsilon}\setminus B_{\varepsilon}} (u_1^m+u_2^m) |\nabla\eta_\varepsilon \cdot \nabla\varphi| \, dx 
    \leq c \varepsilon^{-1} \int_{\varepsilon}^{2\varepsilon} r \, dr \to 0 \quad \text{ as } \quad \varepsilon \to 0.
\end{split}
\end{equation*}
By~\eqref{eq:v_estimates}, $|\Delta \eta_\varepsilon| \leq c_0 \varepsilon^{-2}$ for some $c_0 >0$, and for the fixed $t>0$ we obtain
\begin{equation*} 
    | K_\varepsilon| \leq \int_{B_{2\varepsilon}\setminus B_{\varepsilon}} \varphi |u_1^m-u_2^m| |\Delta\eta_\varepsilon| \, dx 
    \leq c \varepsilon^{-2} b(t) \int_{\varepsilon}^{2\varepsilon} r^{n-1-\lambda} \, dr \to 0 \quad \text{ as } \quad \varepsilon \to 0.
\end{equation*}
By Lemma~2.1 and~\ref{eq:xi}, it also holds that
$\lvert \partial_\tau\eta_\varepsilon\rvert \le c_0 \varepsilon^{-1/\alpha}$ with $\alpha=1$, hence,
\begin{align*}
    |L_\varepsilon| \leq \int_{B_{2\varepsilon}\setminus B_\varepsilon}
    \lvert \partial_\tau \eta_\varepsilon\rvert\,
    \varphi\,\lvert u_1-u_2\rvert\,dx
    &\le
    c\,\varepsilon^{-1}
    \int_{\varepsilon}^{2\varepsilon} r^{n-1-\frac{n-2}{m}}\,dr \\
    &\le
    c\,\varepsilon^{\,\frac{(n-1)}{m}(m-m_*)} \; \to \;0
    \qquad \text{as } \varepsilon \to 0,
    \end{align*}
since $m>m_*=(n-2)/(n-1)$.
Finally, by~\eqref{eq:v_estimates}, $|\eta_\varepsilon| \leq 1$, and $m>m_c$, we have
\begin{equation*} 
\begin{split}
    \left| H_\varepsilon \right| 
    &\leq \sup_{B_{2\varepsilon}\setminus B_{\varepsilon}} \varphi \int_{B_{2\varepsilon}\setminus B_\varepsilon} |u_1-u_2| \, dx 
    \leq \sup_{B_{2\varepsilon}\setminus B_{\varepsilon}} \varphi \int_{B_{2\varepsilon}\setminus B_\varepsilon} (u_1+u_2) \, dx \\
    &\leq c \int_{\varepsilon}^{2\varepsilon} r^{n-1-\frac{n-2}{m}} 
    \leq c \, \varepsilon^{\frac{n}{m} (m-m_c)} \to 0 \quad \text{ as } \quad \varepsilon \to 0.
\end{split}
\end{equation*}
Hence, for any nonnegative function $\varphi\in C_0^\infty (\R^n)$ it holds that
\begin{equation}\label{eq:ineq_m<1_main}
   \int_{\R^n} \varphi |u_1-u_2| \, dx 
    \leq \int_0^t \int_{\R^n} |u_1^m-u_2^m| \Delta \varphi \, dx \, d\tau.
\end{equation}
Furthermore, from~\eqref{eq:ineq_m<1_main} we can derive a useful estimate that will be needed later.
In order to do so, we recall the reverse triangle inequality $|a|a|^{m-1}-b|b|^{m-1}|\leq 2 |a-b|^m$ with exponent $m<1$ and $a, b\in\R$.
Together with the H\"older inequality, we obtain
\begin{equation*} 
\begin{split}
    \int_{\R^n} |u_1^m-u_2^m| \Delta \varphi \, dx
    &\leq 2\int_{\R^n} (|u_1-u_2|\varphi)^m |\Delta \varphi|\varphi^{-m} \, dx \\
    &\leq 2 C[\varphi]^{1-m} \left(\int_{\R^n} \varphi |u_1-u_2| \, dx \right)^m, 
\end{split}
\end{equation*}
where
\begin{equation} \label{eq:c_varphi}
    C[\varphi]:= 
    \int_{\R^n} |\Delta \varphi|^{\frac{1}{1-m}} \varphi^{-\frac{m}{1-m}} \, dx.
\end{equation}
Equation~\eqref{eq:ineq_m<1_main} can be now written as 
$$f'(t) \leq 2 C[\varphi]^{1-m} f^m(t)$$ with $f(0)=0$, and so
\begin{equation} \label{eq:ineq_m<1}
    \int_{\R^n} \varphi |u_1-u_2| \, dx 
    \leq  C[\varphi] \left(2(1-m)t\right)^{\frac{1}{1-m}}.
\end{equation}

\textit{Step 3.}
The rest of the proof is the same as the latter part of the proof of Theorem~2.2 by Takahashi and Yamamoto in~\cite{TY21} and Theorem~2.3 by Herrero and Pierre in~\cite{HP}. We present it for completeness.
Set
\begin{equation*}
    w(x,t):= \int_0^t |u_1^m-u_2^m| \, d\tau.
\end{equation*}
Since $u \in C([0,\infty);L^1_{loc}(\R^n))$ 
and $\varphi\in C_0^\infty (\R^n)$, Fubini’s theorem gives
\begin{equation*}
    \int_{\R^n} \varphi |u_1-u_2| \, dx 
    \leq \int_{\R^n} w(x,t) \Delta \varphi(x) \, dx.
\end{equation*}
Then, $\int_{\R^n} w(x,t) \Delta \varphi(x) \, dx \geq 0$ and so $-\Delta w(x,t) \leq 0$ in $\mathcal{D}'(\R^n)$. 
Hence, the following mean value inequality for subharmonic functions holds
\begin{equation*}
    w(z,t) \leq \frac{1}{|B_1| R^n} \int_{B_R(z)} w(x,t) \, dx =: M_R,
\end{equation*}
where $z\in\R^n$, $|B_1|$ is the volume of a unit ball, and $R>0$.
Thus, $u_1 \equiv u_2$ will be proved once we show $M_R \to 0$ as $R \to \infty$.
For $R \geq 1$ we define $\phi_R \in C_0^\infty (\R^n)$ such that $0 \leq \phi_R\leq 1$, $\phi_R=0$ if $|x-z| \geq 2R$, and $\phi_R=1$ if $|x-z| \leq R$. Let $\tilde \phi_R:= \phi_R^k$ for $k>2/(1-m)$. 
We proceed by using the reverse triangle inequality, H\"older inequality, and~\eqref{eq:ineq_m<1} with $C[\tilde \phi_R]$, which was defined in~\eqref{eq:c_varphi}. We obtain
\begin{equation*}
\begin{split}
    M_R &\leq \frac{2}{|B_1| R^n} \int_0^t \int_{B_R(z)} |u_1-u_2|^m \, dx \, d\tau \\
    &\leq \frac{2}{|B_1|^m R^{nm}} \int_0^t \left(\int_{B_R(z)} |u_1-u_2| \, dx \right)^m\, d\tau \\
    &\leq 2 |B_1|^{-m}R^{-nm} \int_0^t \left(\int_{\R^n} \tilde \phi_R |u_1-u_2| \, dx \right)^m\, d\tau \\
    &\leq (2(1-m))^{\frac{1}{1-m}}|B_1|^{-m} R^{-nm} C[\tilde \phi_R]^m t^{\frac{1}{1-m}}.
\end{split}
\end{equation*}
Substituting $x-z=R(y-z)$, it holds that
\begin{equation*}
    C[\tilde\phi_R] = \int_{B_{2R}(z)} |\Delta \tilde\phi_R|^{\frac{1}{1-m}} \tilde\phi_R^{-\frac{m}{1-m}} \, dx 
    = R^{n-\frac{2}{1-m}} \int_{B_{2}(z)} |\Delta \tilde\phi_1|^{\frac{1}{1-m}} \tilde\phi_1^{-\frac{m}{1-m}} \, dy 
    = R^{n-\frac{2}{1-m}} C[\tilde\phi_1].
\end{equation*}
Since $k>2/(1-m)$, we have
\begin{equation*} \label{eq:c_phi_1}
    C[\tilde\phi_1] = \int_{\R^n} |k(k-1) \phi_1^{k(1-m)-2}|\nabla \phi_1|^2 + k \phi_1^{k(1-m)-1} \Delta \phi_1 |^{\frac{1}{1-m}} \, dx < \infty.
\end{equation*}
Thus,
\begin{equation*}
    M_R \leq (2(1-m))^{\frac{1}{1-m}} |B_1|^{-m} R^{-\frac{2m}{1-m}} C[\tilde\phi_1]^m t^{\frac{1}{1-m}} \to 0 \quad \text{ as } R\to \infty.
\end{equation*}
This shows that $u_1 \equiv u_2$, which completes the proof.
\end{proof}



\begin{thebibliography}{99}




\bibitem{Bookholt}
F. D. Bookholt, H. N. Monsuur, S. Gibbs and F. J. Vermolen,
\textit{Mathematical modelling of angiogenesis using continuous cell-based models},
Biomech. Model. Mechanobiol. \textbf{15} (2016), 1577--1600.

\bibitem{CZ}
X. Chen and W. Zhu,
\textit{A mathematical model of regenerative axon growing along glial scar after spinal cord injury},
Computational and Mathematical Methods in Medicine {\bf 2016} (2016), Art. ID 3030454, 9 pp.

\bibitem{DK} 
\textsc{P. Daskalopoulos and C. E. Kenig},
Degenerate Diffusions: Initial Value Problems and Local Regularity Theory,
EMS tracts in mathematics, European Mathematical Society, 2007.

\bibitem{FTY} 
\textsc{M. Fila, J. Takahashi and E. Yanagida},  
Solutions with moving singularities for equations of porous medium type,
Nonlinear Analysis 179 (2019), 237--253.

\bibitem{FMTY} 
\textsc{M. Fila, P. Macková, J. Takahashi and E. Yanagida},  
Moving singularities for nonlinear diffusion equations in two space dimensions,
Journal of Elliptic and Parabolic Equations 6 (2020), 155--169.

\bibitem{FMTY21} 
\textsc{M. Fila, P. Macková, J. Takahashi and E. Yanagida},  
Anisotropic and isotropic persistent singularities of solutions of the fast diffusion equation,
Differential and Integral Equations 35 (2022), 729--748.

\bibitem{FMY2016} 
\textsc{M. Fila, H. Matano and E. Yanagida},
Non-uniqueness of Solutions of a Semilinear Heat Equation with Singular Initial Data,
In Patterns of dynamics, Springer Proc. Math. Stat. 205 (2017), 138--148.

\bibitem{HP}
\textsc{M. A. Herrero and M. Pierre}, 
The Cauchy problem for $u_t=\Delta u^m$ when $0<m<1$,
Trans. Amer. Math. Soc. 291 (1985), 145--158.

\bibitem{Hui2020}
\textsc{K. Hui}, 
Uniqueness and time oscillating behaviour of finite points blow-up solutions of the fast diffusion equation, 
Proceedings of the Royal Society of Edinburgh: Section A Mathematics 150 (2020), 2849--2870.

\bibitem{Lu1}
\textsc{T. Lukkari}, 
The porous medium equation with measure data, 
J. Evol. Equations 10 (2010), 711--729. 

\bibitem{Lu2}
\textsc{T. Lukkari}, 
The fast diffusion equation with measure data, 
NoDEA Nonlinear Differential Equations Appl. 19 (2012), 329--343.

\bibitem{Oleinik1957}
\textsc{O. A. Oleinik}, 
On equations of the unsteady filtration type, 
Dokl. Akad. Nauk SSSR 113:6 (1957), 1210--1213.

\bibitem{P}
\textsc{M. Pierre}, 
Uniqueness of the solutions of $u_t - \Delta \varphi (u) = 0$ with initial datum a measure,
Nonlinear Analysis: Theory, Methods \& Applications 6 (1982), 175--187.

\bibitem{SY2009}
\textsc{S. Sato and E. Yanagida},
Solutions with moving singularities for a semilinear parabolic equation, 
J. Differential Equations 246 (2009), no. 2, 724--748.

\bibitem{TY21}
\textsc{J. Takahashi and H. Yamamoto}, 
Infinite-time incompleteness of noncompact Yamabe flow,
Calculus of Variations and Partial Differential Equations 61 (2022), no.6, Paper No. 212, 24 pp.

\bibitem{TY15}
\textsc{J. Takahashi and E. Yanagida}, 
Time-dependent singularities in the heat equation, 
Commun. Contemp. Math. 14 (2015), no. 3, 969--979.

\bibitem{vazquez1} 
\textsc{J. L. V{\'a}zquez},
The Porous Medium Equation: Mathematical Theory,
Oxford Mathematical Monographs, Clarendon Press, 2006.


\bibitem{vazquez2} 
\textsc{V{\'a}zquez, J. L.},
Smoothing and Decay Estimates for Nonlinear Diffusion Equations Equations of Porous Medium Type,
Oxford Lecture Notes in Mathematics and Its Applications 33, Oxford University Press, 2006.


\end{thebibliography}
\end{document}